\documentclass[11pt]{amsart} 

\usepackage[english]{babel}
\usepackage[utf8]{inputenc}
\usepackage[T1]{fontenc}
\usepackage{comment}
\usepackage{amsmath}
\usepackage{amssymb}
\usepackage{amsfonts}
\usepackage{amssymb}
\usepackage{amsthm}
\usepackage{amscd}
\usepackage{xcolor}
\usepackage{slashed}
\usepackage[hidelinks]{hyperref}
\usepackage{geometry}
\usepackage{graphicx}
\usepackage{float}

\newcommand{\setsymbol}[1]{\ensuremath{\mathbb{#1}}}%
\newcommand{\R}{\setsymbol{R}}%

\newcounter{Counter}
\newtheorem{definition}{Definition}[section]

\newtheorem{theorem}[definition]{Theorem}
\newtheorem{proposition}[definition]{Proposition}

\newtheorem{remark}[definition]{Remark}

\title{Stability of Llarull's theorem in all dimensions}

\author{Sven Hirsch}
\address{Institute for Advanced Study, 1 Einstein Drive, Princeton, NJ 08540, USA}
\email{sven.hirsch@ias.edu}
\author{Yiyue Zhang}
\address{Department of Mathematics, University of California, Irvine, CA, 92697, USA}
\email{yiyuez4@uci.edu}

\begin{document}

\maketitle

\begin{abstract}
Llarull's theorem characterizes the round sphere $S^n$ among all spin manifolds whose scalar curvature is bounded from below by $n(n-1)$.
In this paper we show that if the scalar curvature is bounded from below by $n(n-1)-\varepsilon$, the underlying manifold is $C^0$-close to a finite number of spheres outside a small bad set.
This completely solves Gromov's spherical stability problem.
\end{abstract}

\section{Introduction}

Recently, there has been much interest in analyzing the stability of various scalar curvature geometry results.
Roughly speaking, these can be divided into two categories:
on the one hand, there are results based on the level-set method such as \cite{ABK2, ABK, DongSong, KKL}.
On the other hand, there are results exploiting the structure of various special settings such as graphs \cite{CGP, HLS} and spherical symmetry \cite{BKS, LeeSormani}.
One particular question is M.~Gromov's \emph{spherical stability problem} \cite[page 20]{Gromov2} concerning the stability of Llarull's theorem \cite{Llarull}.
In this article, we introduce another method to analyze stability: using spin geometry, we confirm M.~Gromov's conjecture in all dimensions.

\begin{theorem}\label{New main theorem}
    Let $(M^n,g)$ be a smooth compact spin manifold without boundary. 
    Let $f:(M^n,g)\to (S^n,g_{0})$ be an area non-increasing smooth map to the round sphere with non-zero degree, and set $N=\max\{\#(f^{-1}(p))|p\in S^n \text{ is a regular value of the map}\}.$
   Suppose that $R_g\ge n(n-1)-\varepsilon$ for some $0\le \varepsilon\le \tfrac 1{CN^2}$ where $C=C(n)$ is an explicit constant. Then we can decompose $M^n$ into smooth sets $\mathcal G,\mathcal B_1,\mathcal B_2$, where 
\begin{align} \label{mainineq}
    |g-f^\ast g_0|_g\le \sqrt \varepsilon
\end{align}
on the good set $\mathcal G$, and where the bad sets $\mathcal B_1,\mathcal B_2$ are small in the sense that
\begin{align}\label{eq:small bad set}
\begin{split}
    |\partial \mathcal B_1|_g+  |\mathcal  B_2|_{g}+|S^n\setminus f(\mathcal G)|_{g_0}\le &CN\sqrt \varepsilon.
    \end{split}
\end{align}
\end{theorem}

\begin{figure}
\includegraphics[totalheight=4cm]{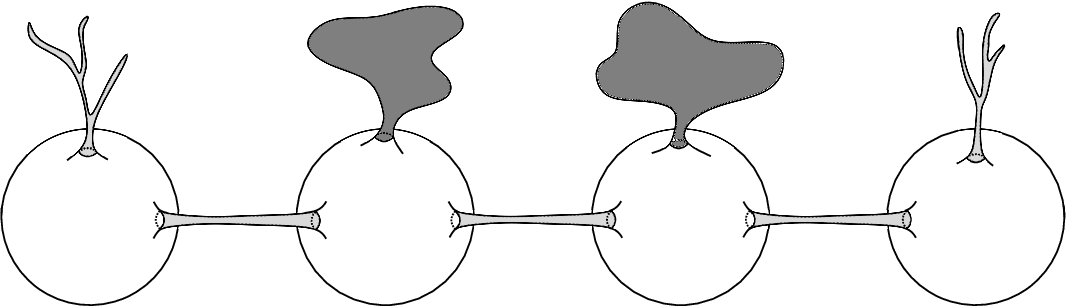}
\caption{
A schematic depiction of the statement of Theorem \ref{New main theorem}. On the good set (white) the manifold is close to a round sphere, and the bad sets $\mathcal B_1$ (dark grey) as well as $\mathcal B_2$ (light grey) are small measured with respect to boundary area and volume respectively.
}
\label{FigureNew}
\end{figure}

The idea of removing a bad set has also been used by C.~Dong and A.~Song \cite{DongSong} for the stability of the positive mass theorem in dimension 3.
They excised a small set (measured by the area of its boundary) and demonstrated that the remainder is Gromov-Hausdorff close to Euclidean space\footnote{We refer to C.~Sormani's survey \cite{Sormani} for a detailed overview of relevant results and conjectures regarding scalar curvature stability.}.

\begin{remark}
It is necessary to have the constant $N$ appearing in equation \eqref{eq:small bad set} since we can concatenate an arbitrary large number of spheres without violating the scalar curvature lower bound.
    In case $f$ does not change orientation in the sense of $\operatorname{det}(Df)\ge0$, we can replace $N$ by $\operatorname{deg}(f)$ in the statement of Theorem \ref{New main theorem}.
\end{remark}

In case $f$ is the identity map, we obtain a slightly stronger result:

\begin{theorem}\label{New main corollary}
Under the same assumptions as in Theorem \ref{New main theorem}, suppose additionally that $f=\operatorname{Id}$. 
Then we can decompose $S^n$ into smooth sets $\mathcal G$ and $\mathcal B$ such that
\begin{align}
    |g-g_0|_g\le \sqrt\varepsilon
\end{align}
on $\mathcal G$, and 
    \begin{equation}
    |\mathcal{B}|_{g_0}\le C\sqrt \varepsilon
    \end{equation}
    for some explicit constant $C=C(n)$.
\end{theorem}

We remark that it is necessary to measure the volume of the bad set with respect to $g_0$ instead of $g$, and the statement is false otherwise.

\subsection{Controlling the geometry of $M$ directly}

In Theorem \ref{New main theorem}, the smallness of the bad set is expressed in a rather subtle way.
This is necessary as seen from Figure \ref{FigureNew} and the examples constructed by P.~Sweeney in \cite{Sweeney}.

If we assume an additional bound on the normalized Sobolev constants\footnote{See Definition \ref{SN} and Remark \ref{remark sobolev}.} $C_{S_\alpha}^\ast$, we obtain more precise control over the bad set.
Intriguingly, these assumptions still allow allow for certain \emph{bags of gold} or \emph{bubbles}, cf. Figure \ref{Figure} below.

\begin{theorem}\label{thm main 3}
Let $\alpha\in[n,\infty]$, $n\ge3$, and $(M^n,g)$ a smooth compact spin manifold without boundary admitting an area non-increasing map onto $(S^n,g_0)$ with non-zero degree.
Let $\varepsilon\in(0,0.1)$ with $\sqrt\varepsilon {C^*_{S_\alpha}}\le 1$, and suppose that 
$\|(R_{g}- n(n-1))_-\|_{L^{\tfrac \alpha2}(M^n)}\le \varepsilon|M^n|^\frac{2}{\alpha}_g$, where $x_-=\max\{-x,0\}$.
Then there exists a decomposition of $M$ into the smooth sets $\mathcal G$ and $\mathcal B$ such that 
\begin{align} \label{g-g0}
      |g-f^\ast g_0|_g\le n\sqrt\varepsilon
\end{align}
on $\mathcal G$, and
\begin{align}
    |\mathcal B|_g\le  \left[257\varepsilon^\frac{\alpha}{2(\alpha+2)}+384({C^*_{S_\alpha}})^\frac{2\alpha}{\alpha-2}\cdot\varepsilon^\frac{\alpha}{\alpha-2}\right] |M^n|_g.
\end{align}
\end{theorem}

\begin{remark}
For $\alpha=\infty$, this yields a stability result involving a pointwise lower bound for the scalar curvature and the Poincar\'e constant.
\end{remark}

\begin{figure}
\includegraphics[totalheight=4cm]{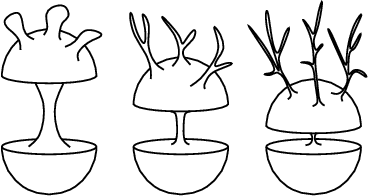}
\caption{
A sequence of manifolds $(S^n,g_i)$ with uniformly bounded Poincar\'e constants. We construct such an example explicitly in Proposition \ref{prop:example}.
}
\label{Figure}
\end{figure}
    
With similar methods we can also obtain intrinsic flat convergence.
This notion of convergence has been introduced by C.~Sormani and S.~Wenger \cite{SormaniWenger} and builds upon the theory of currents in metric spaces developed by L.~Ambrosio and B.~Kirchheim \cite{AmbrosioKirchheim}.
The popularity of intrinsic flat convergence stems from its ability to control \emph{splines} or \emph{trees} which can appear in the limiting process as shown in Figure \ref{Figure}.
These splines make it impossible to expect $C^0$ or even Gromov-Hausdorff convergence on the entire manifold.

\begin{theorem}\label{thm main 2}
Let $\alpha\in[n,\infty]$, and let $g_i$ be a sequence of smooth metrics on $S^n$, $n\ge3$, such that
    \begin{enumerate}
        \item $\|(R_{g_i}- n(n-1))_-\|_{L^{\tfrac \alpha2}(S^n)}\le \tfrac1i |S^n|_{g_i}^\frac{2}{\alpha}$,
        \item $g_i\ge g_0$,
        \item the normalized Sobolev constants\ $C_{S_\alpha}^\ast$ of $(S^n,g_i)$ are uniformly bounded from above,
        \item the diameters of $(S^n,g_i)$ are uniformly bounded from above.
    \end{enumerate}
Then $(S^n,g_i)$ converges to $(S^n,g_0)$ in the intrinsic flat sense.
\end{theorem}

Assuming a Cheeger constant bound instead of a Poincar\'e or Sobolev constant bound, a similiar version of this statement has already been established in dimension 3 by B.~Allen, E.~Bryden and D.~Kazaras in \cite{ABK}.
Their proof relies on an integral formula for spacetime harmonic functions \cite{HKKZ} due to D.~Kazaras, M.~Khuri and the authors which originates from the proof of the spacetime positive mass theorem \cite{HKK, HirschZhang}.

\textbf{Acknowledgements:} Part of this work was carried out while the authors attended a conference at the Simons Center for Geometry and Physics and the authors are grateful for the stimulating research environment. 
SH was supported by the National Science Foundation under Grant No. DMS-1926686, and by the IAS School of Mathematics. 
The authors also thank Brian Allen, Edward Bryden, Simone Cecchini, Bernhard Hanke, Demetre Kazaras, Richard Schoen, Christina Sormani, and Rudolf Zeidler for insightful discussions and their interest in this work.


\section{Cheeger, Sobolev, and Poincar\'e constants}\label{SS:Cheeger}

In this section, we recall the definitions of Cheeger, Sobolev, and Poincar\'e constants and discuss their interrelationships.
Throughout this section, $(M^n,g)$ is a smooth compact Riemannian manifold without boundary.
\subsection{Definitions and preliminaries}
\begin{definition}
    The Poincar\'e constant of $(M^n,g)$ is defined by
    \begin{align}
        C_P=\sup_{u\in W^{1,2}(M)}\frac{\inf_{a\in \R}\| u-a\|_2^2}{\|\nabla u\|_2^2}.
    \end{align}
\end{definition}
We remark that $C_P^{-1}$ is the first eigenvalue of the Laplace operator on $M^n$.

Before we define the Sobolev constants $C_{S_\alpha}$, we first introduce two other, more commonly used geometric constants, see for instance \cite[Definition 2.4]{DWZ} by X.~Dai, G.~Wei and Z.~Zhang, and \cite{P.Li} by P.~Li.

\begin{definition} 
For $\alpha\in[n,\infty]$, the Neumann $\alpha$-isoperimetric constant $\operatorname{IN}_\alpha$ of $M$ is defined by 
    \begin{equation}\label{Cheeger def}
        \operatorname{IN}_\alpha(M)=\sup_\Sigma \frac{\min\{|\Omega_1|,|\Omega_2|\}^{1-\frac{1}{\alpha}}}{|\Sigma|},
    \end{equation}
    where $\Sigma$ is any smooth hypersurface dividing $M$ into the regions $\Omega_1$ and $\Omega_2$.
The Neumann $\alpha-$Sobolev constant of $M$ is defined by 
\begin{equation}\label{Sobolev constant SN}
    \operatorname{SN}_{\alpha}=\sup_{f\in W^{1,1}(M)}\frac{\inf_{a\in \R}\|f-a\|_{\frac{\alpha}{\alpha-1}}}{\|\nabla f\|_1}.
\end{equation}
\end{definition}

The Neumann $\infty$-isoperimetric $\operatorname{IN}_\infty$ is usually referred to as the Cheeger constant.
We remark that some authors refer to $\operatorname{IN}_\infty^{-1}$ as the Cheeger constant instead.

Next, we recall from \cite[Theorem 9.2]{P.Li}. 
\begin{proposition}
    For all $n\le \alpha\le \infty$, we have 
\begin{equation}
    \frac{1}{2}\operatorname{IN}_{\alpha}(M)\le \operatorname{SN}_{\alpha}(M)\le \operatorname{IN}_{\alpha}(M).
\end{equation}
\end{proposition}

\begin{definition}\label{SN}
    For $\alpha\in[n,\infty]$, the Sobolev constant $C_{S_\alpha}$ and the normalized Sobolev constant $C^*_{S_\alpha}$  are defined by
        \begin{equation}\label{Sobolev constant}
    C_{S_\alpha}=\sup_{f\in W^{1,2}}\frac{\inf_{a\in \R}\|f-a\|_{\frac{2\alpha}{\alpha-2}}}{\|\nabla f\|_2}\quad\text{and}\quad C_{S_\alpha}^*=C_{S_\alpha}|M|^\frac{1}{\alpha}.
\end{equation}
\end{definition}

In Section \ref{SS:Cheeger comparison},  we will delve deeper into the case $\alpha=\infty$ and compare the Cheeger constant to the Poincar\'e constant.
We will observe that a bound on the former implies a bound on the latter.
Similarly, a bound on $\operatorname{SN}_{\alpha}$ implies a bound on $C_{S_\alpha}$.
For the reader's convenience, we sketch the following proof of R.~Schoen \cite{Schoen}.

\begin{proposition} 
    Let $p\in [1,\alpha)$.
    For any $f\in C^\infty(M)$, we have
    \begin{equation}
       \inf_{a\in\R} \|f-a\|_\frac{p\alpha}{\alpha-p}\le \frac{p(\alpha-1)}{\alpha-p}\operatorname{IN}_\alpha\|\nabla f\|_p.
    \end{equation}
\end{proposition}
\begin{proof}
    Let $\bar{f}$ be the constant such that $\{f\ge \bar{f}\}\ge \frac{1}{2}|M|$ and
    $\{f\le \bar{f}\}\ge \frac{1}{2}|M|$. 
    Then we have
    \begin{equation} \label{INalpha}
        \operatorname{IN}_{\alpha}(M)=\sup_{f\in C^\infty(M)}\frac{\|f-\bar{f}\|_{\frac{\alpha}{\alpha-1}}}{\|\nabla f\|_1}.
    \end{equation}
    Let $v=(f-\bar{f})_+$ and observe that $\{v\ge 0\}\ge \frac{1}{2}|M|$ and $\{v\le 0\}\ge \frac{1}{2}|M|$.
    Hence, we may apply \eqref{INalpha} with $\bar{v}=0$.
    Combining this with H\"older's inequality yields
    \begin{equation} \label{sIN}
        \left(\int_M v^\frac{s\alpha}{\alpha-1}dV\right)^\frac{\alpha-1}{\alpha}
        \le s\operatorname{IN}_\alpha\int_M v^{s-1}|\nabla v|dV
        \le s\operatorname{IN}_{\alpha}
        \left(\int_M v^{(s-1)\frac{p}{p-1}}dV\right)^\frac{p-1}{p}\|\nabla v\|_{p}
    \end{equation}
where we set $s=\frac{p(\alpha-1)}{\alpha-p}$ which ensures $\frac{s\alpha}{\alpha-1}=\frac{p(s-1)}{p-1}$.
Similarly, we obtain
    \begin{equation}
        \left(\int_M (f-\bar{f})_-^\frac{s\alpha}{\alpha-1}dV\right)^\frac{\alpha-1}{\alpha}
        \le s\operatorname{IN}_{\alpha}
        \left(\int_M (f-\bar{f})_-)^{(s-1)\frac{p}{p-1}}dV\right)^\frac{p-1}{p}\|\nabla (f-\bar{f})_-\|_{p}.
    \end{equation} 
Adding the above two inequalities gives
    \begin{equation}
       \left( \int_M |f-\bar{f}|^{\frac{p\alpha}{\alpha-p}}dV\right)^{\frac{1}{p}-\frac{1}{\alpha}}\le \frac{p(\alpha-1)}{\alpha-p}\operatorname{IN}_{\alpha}\|\nabla f\|_p 
    \end{equation}
    which finishes the proof.
\end{proof}

\begin{remark}\label{remark sobolev} Since $g_i\ge g_0$ implies a lower volume bound, Theorem \ref{thm main 2} still holds with an upper bound on $C_{S_\alpha}$ instead of $C_{S_\alpha}^\ast$.
Moreover, according to the above proposition, Theorem \ref{thm main 2} also holds with an upper bound on $\operatorname{SN}_\alpha$ and $\operatorname{IN}_\alpha$. 
\end{remark}

\subsection{Comparison between the Cheeger and the Poincar\'e constant}\label{SS:Cheeger comparison}

As mentioned above, a uniform upper bound on the Cheeger constants implies a uniform upper bound on the Poincar\'e constants.
This is known as Cheeger's inequality \cite{Cheeger}.

According to Buser's inequality \cite{Buser2}, the reverse inequality also holds as long as a uniform lower Ricci curvature bound is imposed.
This condition is necessary as demonstrated in \cite[Section 4]{Buser1} by P.~Buser where the conformal invariance of the Dirichlet integral in dimension 2 is exploited.
In the proposition below, we adapt this example to higher dimensions.
In particular, a Poincar\'e constant upper bound does not prevent the presence of bags of gold, cf. Figure \ref{Figure}.

\begin{proposition}\label{prop:example}
There exists a family of metrics $g_\delta $ on $S^n$ such that the Poincar\'e constants of $g_\delta$ are uniformly bounded from above while the Cheeger constants diverge to $\infty$.
\end{proposition}

\begin{proof}
Let $g_0$ be the round metric, fix an equator $S^{n-1}\subset S^n$, and let $\rho$ be the signed distance function (w.r.t. $g_0$) to this equator.
For each $\delta\in(0,\tfrac12]$, we define $g_\delta=ug_0$ where $u=u_\delta(\rho)$ is a function depending only on $\rho$. 
More precisely, we prescribe that 
\begin{itemize}
    \item $u(\rho)=u(-\rho)$,
    \item $u(\rho)=1$ for $\rho\in [\delta,\tfrac\pi2]$,
    \item $u$ attains its minimum at $\rho=0$ with $u(0)=\delta^{\tfrac{1}{n-1}}$,
    \item $u\le1$ for all $\rho\in[-\tfrac\pi2,\tfrac\pi2]$.
\end{itemize}
It is easy to see that the Cheeger constants of $(S^n,g_\delta)$ diverge to $\infty$ for $\delta\to0$ in view of Definition \eqref{Cheeger def}.

To see that the Poincar\'e constants stay bounded for $\delta\to0$, we proceed as follows.
Let $\psi$ be the first eigenfunction of the Laplacian and recall that it realizes the Poincar\'e constant, i.e. $C_P=\tfrac{\int_{S^n} \psi^2dV_g}{\int_{S^n}|\nabla\psi|^2dV_g}$.
Here we omit any $\delta$ subscripts to declutter the notation.
First, observe that $\psi$ inherits the symmetries of $u$, i.e. $\psi$ depends only on $\rho$, and $\psi(\rho)=\psi(-\rho)$. 
Moreover, we may scale $\psi$ such that $\psi(\pm\frac{\pi}{2})=\pm1$.
Next, we claim that $\psi$ is monotonically increasing in $\rho$ which in particular implies that $|\psi|\le1$.
Suppose $\psi$ is not monotone.
Then there exist various local minima and maxima which will be attained at $-a_{k},-a_{k-1},\dots,-a_1,a_1,\dots,a_{k-1},a_k$.
We denote with $b_k$ the function values $b_k=\psi (a_k)$. 
Note that $b_k=-b_{-k}$.
Let $b_l$ be the maximal value of $\{b_1,\dots, b_k\}$.
We distinguish two cases.
First suppose that $b_l\ge1$.
In this case, consider the function $\overline \psi$ defined by $\overline \psi(x)=b_k$ for $x\ge a_k$, and $\overline \psi(x)=\psi(x)$ for $x<a_k$.
Then $\int_{S^n}\overline \psi^2dV_g>\int_{S^n}\psi^2dV_g$ and $\int_{S^n}|\nabla \overline\psi|^2dV_g<\int_{S^n}|\nabla\psi|^2dV_g$.
This is a contradiction to the assumption that $\psi$ realizes the Poincar\'e constant.
Suppose now that that $b_l<1$. 
In this case, we can find a point $c\in (a_l,1)$ such that $\psi(c)=b_l$ and $\psi(x)<b_l$ for $x\in(a_l,c)$.
Now we proceed as above by setting $\overline\psi(x)=b_l$ for $x\in [a_l,c]$ which again leads to contradiction.
Hence, $\psi$ is monotone and we have $|\psi|\le 1$.

To estimate the Poincar\'e constants of $(S^n,g_\delta)$ we need to make another case decomposition.
First, suppose that $\psi(\delta)\le\frac12$.
To show that the Poincar\'e constant is uniformly bounded (in $\delta$), we need to estimate $\int_{S^n}\psi^2dV_g$ from above and $\int_{S^n}|\nabla\psi|^2dV_g$ from below. 
Clearly, $\int_{S^n}\psi^2dV_g\le |S^n|_{g_0}$ since $|\psi|\le1$.
Define $\overline \psi$ by setting $\psi(\rho)=0$ for $\rho\in[-\delta,\delta]$, and $\overline\psi(\rho)=a\psi(\rho)-b$ otherwise.
Here, the constants $a\in[1,2]$ and $b\in[0,1]$ are chosen such that $\overline\psi$ is Lipschitz and $\overline \psi(\pm1)=\pm1$.
Then $\overline \psi$ is a valid competitor for computing the Poincar\'e constant for $(S^n,g_0)$.
Moreover, $\int_{S^n}|\nabla \psi|^2dV_g\ge\tfrac{1}{4}\int_{S^n}|\nabla\overline\psi|^2dV_{g_0}$.
Thus, we have obtained a uniform bound for $C_P$ in this case.

Next, suppose that $\psi(\delta)>\frac12$.
Again, we obtain the bound $\int_{S^n}\psi^2dV_g\le |S^n|_{g_0}$ and it remains to estimate $\int_{S^n}|\nabla\psi|^2dV_g$.
We compute
\begin{equation}
\int_{S^n}|\nabla\psi|^2dV_g\ge\int_{-\delta\le\rho\le\delta}|\nabla \psi|^2dV_g\ge \int_{-\delta}^\delta |S^{n-1}|_{g_0} \left(\sqrt{1-\delta^2}u(0)\right)^{n-1}|\nabla \psi|^2d\rho\ge C(n) \delta \int_{-\delta}^\delta |\nabla \psi|^2d\rho 
\end{equation}
where $C(n)$ is a constant depending only on $n$.
Since $u\le1$, we have $|\nabla\psi|^2\ge (\partial_\rho\psi)^2$.
Hence,
\begin{align}
    \int_{-\delta}^\delta|\nabla \psi|^2d\rho\ge\int_{-\delta}^\delta (\partial_\rho\psi)^2d\rho\ge \frac1{2\delta}\left(\int_{-\delta}^\delta \partial_\rho\psi d\rho\right)^2> \frac1{2\delta},
\end{align}
and the result follows. 
\end{proof}

We remark that $(S^n,g_\delta)$ does not converge to $(S^n,g_0)$ in the intrinsic flat convergence, cf. \cite[Example 2.5]{APS}. 
We also note that $g_\delta\ge g_0$ does not hold, and the above example is therefore not a counter example to Theorem \ref{thm main 2}.

A similar example can be constructed where the Neumann $\alpha$-Sobolev constants $\operatorname{SN}_\alpha$ go to infinity, but the Sobolev constants $C_{S_\alpha}$ remain bounded. 
The key is to make the cylinder connecting the two hemispheres very short so the gradient of $\phi$ cannot accumulate in this region.


\section{Proof of Theorem \ref{New main theorem} and Theorem \ref{New main corollary}}

\begin{proof}[Proof of Theorem \ref{New main theorem}] 
To simplify the notation we abbreviate $|\cdot|_g$ by $|\cdot|$ below.
We define the function $F:S^n\to \mathbb R$ via
\begin{equation}
    F(p)=\sup_{x\in M^n}\{|\phi|^2(x)| f(x)=p\}
\end{equation}
Since the critical values of $f$ on $S^n$ are a set of measure zero by Sard's theorem, $F$ is well defined almost everywhere on
$S^n$. Moreover, $F$ is locally Lipschitz on the set of regular values of $f$, and $F$ is $W^{1,1}$.

Let $0<\lambda_1\le \lambda_2\le\cdots \le \lambda_n$ be the square roots of eigenvalues of $g$ under the pullback metric $f^*g_0$. 
In case, the differential $df$ degenerates, we take the inverses of the eigenvalues of $f^*{g_0}$ with respect to $g$. 
In particular, we obtain $\lambda_n=\infty$ in this setting.
The assumption that $f$ is area non-increasing implies $\lambda_i\lambda_j\ge 1$, for any $1\le i<j\le n$.

Next, we define the set
\begin{equation}
    \mathcal{B}_0=\{x\in M^n| \lambda_{n-1}\lambda_{n}\ge 1+\tfrac1{4n}\sqrt{\varepsilon}\}
\end{equation}
and set $\mathcal{G}_0=\mathcal{B}_0^c$. 
We note that if $x\in \mathcal G_0$, we have $1-\tfrac1{4n}\sqrt\varepsilon<\lambda_i <1+\tfrac1{4n}\sqrt\varepsilon$ for all $1\le i\le n$ since $n\ge3$. In particular $|g-f^\ast g_0|_g\le \sqrt\varepsilon$ on $\mathcal G_0$ for $\varepsilon$ sufficiently small. 
We estimate on $\mathcal{B}_0$
\begin{equation} \label{BRg}
\sum_{j\neq l}\frac{1}{\lambda_l\lambda_j}-R_g\le \frac{2}{\lambda_{n-1}\lambda_n}-2+\varepsilon\le 2(1+\tfrac1{4n}\sqrt\varepsilon)^{-1}-2+\varepsilon=-\tfrac1{2n}\sqrt\varepsilon (1+\tfrac1{4n}\sqrt{\varepsilon})^{-1}+\varepsilon\le -\tfrac1{4n}\sqrt\varepsilon
\end{equation}
for $\varepsilon$ sufficiently small depending only on $n$.
Next, we use Llarull's integral formula, cf. Appendix \ref{appendix:llarull}, which in even dimensions states
\begin{equation}
    \int_{M^n} |\phi|^2\left[\sum_{j\neq l}\frac{1}{\lambda_j\lambda_l}-R_g\right]dV_g\ge 4\int_{M^n} |\nabla \phi|^2dV_g.
\end{equation}
For the adjustments needed in the odd dimensional case see Appendix \ref{SS:odd}.
In our setting, Llarull's formula implies
\begin{equation} \label{G0B0}
\int_{\mathcal{G}_0}\varepsilon |\phi|^2 dV_g\ge \int_{\mathcal{B}_0}\frac1{4n} \sqrt\varepsilon|\phi|^2dV_g+4\int_{M^n} |\nabla |\phi||^2dV_g.
\end{equation}
Hence, we have for sufficiently small $\varepsilon$
\begin{equation} \label{GS1}
   2\varepsilon \int_{\mathcal{G}_0}|\phi|^2 dV_g\ge \int_{M^n} \left(\varepsilon|\phi|^2+4|\nabla |\phi||^2\right)dV_g.
\end{equation}
Let $\mathring{\nabla}$ be the connection on $(S^n,g_0)$.
Choosing $\varepsilon>0$ sufficiently small, depending only on $n$, we may impose that $\lambda_1\cdots\lambda_n\le2$ holds on $\mathcal{G}_0$.
Moreover, we have $\lambda_1\cdots\lambda_{n-1}\ge1$ on the entire manifold $M^n$ since $f$ is area non-increasing.
Thus,
\begin{equation} \label{fGS}
    \begin{split}
        4N\varepsilon\int_{f(\mathcal{G}_0)} FdV_{g_0}\ge&  2\varepsilon \int_{\mathcal{G}_0}|\phi|^2 dV_g 
        \\ \ge& \int_{M^n} \left(\varepsilon|\phi|^2+4|\nabla |\phi||^2\right)dV_g
       \\ \ge& \int_{S^n} \left(\varepsilon F\lambda_1\cdots\lambda_n+4|\mathring{\nabla} F^\frac{1}{2}|_{g_0}^2\lambda_1\cdots\lambda_{n-1}\lambda_n^{-1}\right)dV_{g_0}
              \\ =& \int_{S^n}\sqrt\varepsilon \lambda_1\cdots \lambda_{n-1}\left[\lambda_n\sqrt\varepsilon F+(\lambda_n\sqrt\varepsilon F)^{-1}|\mathring{\nabla} F|_{g_0}^2\right]dV_{g_0}
        \\ \ge& 
       \int_{S^n} \sqrt\varepsilon |\mathring{\nabla}F|_{g_0}dV_{g_0}
       \\ \ge & C_0\sqrt \varepsilon \inf_{a\in\R}\int_{S^n} \left|F-a\right|dV_{g_{0}}.
    \end{split}
\end{equation}
Here, the first inequality follows from the definition of $F$, the second inequality follows from equation \eqref{GS1}, and in the third we are replacing $g_0$ by $g$.
Moreover, we have used Young's inequality and the Poincar\'e inequality where we have denoted with $C_0=C_0(n)$ the Cheeger constant of $(S^n,g_0)$.
By \eqref{fGS}, $F$ is almost constant on $S^n$, and after rescaling $\phi$, we may choose $a=1$. 
Moreover, estimating $F\le |F-1|+1$, we obtain
\begin{equation}
    4N\varepsilon\int_{S^n} (|F-1|+1)d V_{g_0}\ge C_0\sqrt\varepsilon\int_{S^n} \left|F-1\right|dV_{g_0}.
\end{equation}
Rearranging gives 
\begin{equation} \label{F12}
    \int_{S^n} \left|F-1\right|dV_{g_0}\le \frac{4N\varepsilon}{C_0\sqrt\varepsilon-4N\varepsilon} |S^n|_{g_0}, \quad \text{and}\quad \int_{S^n}  F dV_{g_0}\le \frac{C_0\sqrt\varepsilon}{C_0\sqrt\varepsilon-4N\varepsilon} |S^n|_{g_0}
\end{equation}
where we impose that $\varepsilon$ is sufficiently small and satisfies $4N\varepsilon\le \frac{1}{3}C_0\sqrt\varepsilon$.
Combining this with the first line of \eqref{fGS} yields
\begin{equation} \label{GVN}
    \int_{\mathcal{G}_0}|\phi|^2 dV_g\le 2N\int_{f(\mathcal{G})} F dV_{g_0}\le \frac{2NC_0\sqrt\varepsilon}{C_0\sqrt\varepsilon-4N\varepsilon} |S^n_{g_0}|\le 3N |S^n|_{g_0}.
\end{equation}
Note that the first inequality in \eqref{F12} implies 
\begin{equation}
    \int_{S^n}|F-1|dV_{g_0}\le 6NC_0^{-1}\sqrt\varepsilon|S^n|_{g_0}.
\end{equation}
Hence,  for any $c\in (0,1)$,
\begin{equation}\label{FCe}
    |\{|F-1|\ge c\}|_{g_0}\le 6Nc^{-1}C_0^{-1}\sqrt\varepsilon|S^n|_{g_0}.
\end{equation}

According to line  \eqref{G0B0}, we have 
\begin{equation} \label{G0phi}
  \varepsilon\int_{\mathcal{G}_0} |\phi|^2 dV_g\ge 4\int_{M^n}|\nabla|\phi||^2 dV_g, \quad \text{and}\quad  \int_{\mathcal{G}_0}|\phi|^2 dV_g \ge\frac{1}{2} \int_{M^n} |\phi|^2 dV_g.
\end{equation}
We denote with $\Sigma_t=\{ |\phi|^2=t\}$ be the level set of $|\phi|^2$ in $M^n$.
Multiplying these two inequalities in \eqref{G0phi} together and applying Cauchy-Schwarz inequality, then using the coarea formula yields
\begin{equation} \label{levelset}
    \sqrt\varepsilon\int_{\mathcal{G}_0}|\phi|^2 dV_g\ge \frac{1}{\sqrt{2}}\int_{M^n} |\nabla |\phi|^2| dV_g= \frac{1}{\sqrt{2}}\int^{\overline{t}}_{\underline{t}} |\Sigma_t| dt,
\end{equation}
where $\overline{t}$ and $\underline{t}$ are the maximal and minimal values of $|\phi|^2$. Combining Line \eqref{GVN} and \eqref{levelset}, we obtain 
\begin{equation}
    \int^{\overline{t}}_{\underline{t}} |\Sigma_t| dV_g \le 3\sqrt{2}N\sqrt\varepsilon |S^n|_{g_0}.
\end{equation}
Consequently, we can find $t_1\in[\frac{5}{4},\frac{3}{2}]$, $t_2\in[\frac{1}{2},\frac{3}{4}]$ such that
\begin{equation} \label{partial B1}
    |\Sigma_{t_1}|\le 12\sqrt{2}N\sqrt\varepsilon|S^n|_{g_0},\quad 
    |\Sigma_{t_2}|\le 12\sqrt{2}N\sqrt\varepsilon|S^n|_{g_0}.
\end{equation}
Moreover, we can impose by Sard's theorem that $\Sigma_{t_1}$ and $\Sigma_{t_2}$ are smooth.

Next, let $\Omega_{[t_1,t_2]}=\{x\in M^n : |\phi|^2\in[t_1,t_2]\}$. Note that 
\begin{equation}
    |f(\Omega)|_{g_0}\ge |\{F\in[t_1,t_2]\}|_{g_0}- |\{F> t_2\}|_{g_0}.
\end{equation}
Since $\Omega_{[\tfrac{3}{4},\tfrac{5}{4}]}\subset\Omega$ and using inequality \eqref{FCe}, we have 
\begin{equation} \label{Omega}
    |f(\Omega)|_{g_0}\ge \left|\left\{F\in\left[\frac{3}{4},\frac{5}{4}\right]\right\}\right|_{g_0}-\left|\left\{F>\frac{5}{4}\right\}\right|_{g_0}\ge (1-48NC_0^{-1}\sqrt\varepsilon)|S^n|_{g_0}.
\end{equation}

Next, we estimate with the help of \eqref{G0B0} and \eqref{GVN}
\begin{equation}
    \int_{\mathcal{B}_0} |\phi|^2 dV_g \le 4n\sqrt\varepsilon\int_{\mathcal{G}_0} |\phi|^2 dV_g\le 12nN\sqrt\varepsilon|S^n|_{g_0}.
\end{equation}
Moreover, we have by the definition of $\Omega$
\begin{equation}
    \int_{\mathcal{B}_0} |\phi|^2 dV_g \ge \int_{\mathcal{B}_0\cap \Omega}|\phi|^2 dV_g \ge \frac{1}{2}|\mathcal{B}_0\cap \Omega|.
\end{equation}
We define the sets $\mathcal{B}_1:=M^n\setminus \Omega$, $\mathcal{B}_2=\mathcal{B}_0\cap\Omega$, and $\mathcal{G}:=\Omega\cap\mathcal{G}_0$.
The above two inequalities imply
\begin{equation} \label{B2}
    |\mathcal{B}_2|<24nN\sqrt\varepsilon|S^n|_{g_0}.
\end{equation}
Hence, Theorem \ref{New main theorem} follows from \eqref{partial B1}, \eqref{Omega} 
 and \eqref{B2}.
\end{proof}

\begin{proof}[Proof of Theorem \ref{New main corollary}]
This follows immediately from Theorem \ref{New main theorem} by setting $\mathcal B=\mathcal B_1\cup \mathcal B_2$, and using that $|S^n\setminus \mathcal G|_{g_0}\le C\sqrt\varepsilon$.
\end{proof}


\section{Proofs of Theorem \ref{thm main 3} and Theorem \ref{thm main 2}}

In this section we show Theorem \ref{thm main 3} and Theorem \ref{thm main 2} in even dimensions.
For the odd dimensional case we refer to Section \ref{SS:odd}, or more precisely to equation \eqref{odd}.

\begin{proof}[Proof of Theorem \ref{thm main 3}]
Again, we abbreviate $|\cdot|_g$ by $|\cdot|$.
 Using H\"older's  inequality and our assumption on $R_g$ yields
    \begin{equation}\label{eq B 1}
        \int_{M^n} |\phi|^2(n(n-1)-R_g)_+dV_g\le \|\phi\|_{{\frac{2\alpha}{\alpha-2}}}^{2}\|(n(n-1)-R_g)_+\|_{\frac{\alpha}{2}}\le \varepsilon \|\phi\|_{{\frac{2\alpha}{\alpha-2}}}^{2}|M^n|^{\frac{2}{\alpha}}.
    \end{equation}
Next,  we apply the Sobolev inequality \eqref{Sobolev constant} to find
\begin{equation}\label{eq: possibility to improve}
       \int_{M^n}|\nabla |\phi||^2 dV_g 
        \ge (C_{S_{\alpha}}^*)^{-2}|M^n|^{\frac{2}{\alpha}} \cdot \| |\phi|-1\|^2_\frac{2\alpha}{\alpha-2}.
\end{equation}
Note that we can choose $a=1$ in the Sobolev inequality \eqref{Sobolev constant} by rescaling $\phi$ appropriately.
Combining equation \eqref{eq B 1}, Llarull's integral formula, cf. Theorem \ref{A.1}, and Kato's inequality, we obtain
\begin{equation}
        \varepsilon \|\phi\|_{{\frac{2\alpha}{\alpha-2}}}^{2}|M^n|^{\frac{2}{\alpha}} \ge 4\int_{M^n}|\nabla \phi|^2dV_g\ge 4(C_{S_{\alpha}}^*)^{-2} |M^n|^{\frac{2}{\alpha}}\cdot \| |\phi|-1\|^2_\frac{2\alpha}{\alpha-2}.
\end{equation}
Consequently,
\begin{align}
         \varepsilon\left( \int_{M^n} |\phi|^{\frac{2\alpha}{\alpha-2}}dV_g\right)^{\frac{\alpha-2}{\alpha}}\ge 
4(C_{S_{\alpha}}^*)^{-2}\| |\phi|-1\|^2_\frac{2\alpha}{\alpha-2}.
\end{align}
Therefore,
\begin{equation}
    \varepsilon^{\frac{\alpha}{\alpha-2}} \int_{M^n} |\phi|^{\frac{2\alpha}{\alpha-2}}dV_g\ge 
   \left(\frac{2}{C^*_{S_\alpha}}\right)^{{\frac{2\alpha}{\alpha-2}}}
    \int_{M^n} \left||\phi|-1\right|^\frac{2\alpha}{\alpha-2}dV_g.
\end{equation}
By the generalized mean inequality, we have
\begin{equation}
    |\phi|^{\frac{2\alpha}{\alpha-2}}=(|\phi|-1+1)^{\frac{2\alpha}{\alpha-2}}\le 2^{\frac{\alpha+2}{\alpha-2}}(\left||\phi|-1\right|^{\frac{2\alpha}{\alpha-2}}+1).
\end{equation}
Thus, using $\varepsilon ({C^*_{S_\alpha}})^2\le 1$, we obtain
\begin{equation} \label{fal}
    \int_{M^n} ||\phi|-1|^\frac{2\alpha}{\alpha-2}dV_g\le \frac{\varepsilon_*^\frac{\alpha}{\alpha-2}}{2-\varepsilon_*^\frac{\alpha}{\alpha-2}}|M^n| \le  \varepsilon_*^\frac{\alpha}{\alpha-2}|M^n|,
\end{equation}
where we set $\varepsilon_*=\varepsilon ({C^*_{S_\alpha}})^2$.
Moreover, 
\begin{equation}
\begin{split}
     \int_{M^n}|\phi|^\frac{2\alpha}{\alpha-2}dV_g
     \le& \int_{M^n} 2^\frac{\alpha+2}{\alpha-2} (||\phi|-1|^\frac{2\alpha}{\alpha-2}+1)dV_g
     \\ \le & 2^\frac{\alpha+2}{\alpha-2}(1+  \varepsilon_*^\frac{\alpha}{\alpha-2}) |M^n|.
\end{split}
\end{equation}
Next, we estimate
\begin{equation} \label{C1al C2al}
 \int_{M^n}|\phi|^\frac{2\alpha}{\alpha-2}dV_g\le 2^6|M^n|    \quad \text{and}\quad \int_{M^n}|\phi|^2 dV_g\le (2^6)^\frac{\alpha-2}{\alpha}|M^n|\le 2^6|M^n|.
\end{equation}
where we used H\"older's inequality, $\varepsilon_*\le 1$, and the fact that $\alpha\ge n\ge 3$ implies $2^{\frac{\alpha+2}{\alpha-2}}\le 2^5$.

Let us denote with $\mathcal E$ the set $\mathcal{E}=\{x\in M^n|n(n-1)-R_g\ge \varepsilon^\frac{\alpha+1}{\alpha+2}\}$.
We find
\begin{equation}
    \varepsilon\ge
    \left[\frac{1}{|M^n|}\int_{M^n}[n(n-1)-R_g]_+^\frac{\alpha}{2}dV_g\right]^\frac{2}{\alpha}\ge \varepsilon^\frac{\alpha+1}{\alpha+2}\left(\frac{|\mathcal{E}|}{|M^n|}\right)^\frac{2}{\alpha}
\end{equation}
which implies 
\begin{align}\label{B estimate}
    |\mathcal{E}|\le \varepsilon^{\frac{\alpha}{2(\alpha+2)}}|M^n|.
\end{align}

Next, let $\mathcal{B}$ be the subset of $M^n$ where $\max_{l\neq j}\{\lambda_l\lambda_j\}\ge1+\sqrt\varepsilon$.
Note that inequality \eqref{g-g0} holds on $\mathcal{B}^c$, since $\lambda_n\in [ (1+\sqrt{\varepsilon})^\frac{1}{2},1+\sqrt{\varepsilon}]$, and 
\begin{equation}
|g-f^*g_0|_g=\left|n-\sum_{i=1}^n \lambda_i^{-2}\right|\le n-\lambda_n^{-2}-(n-1)(1+\sqrt{\varepsilon})^{-2}\lambda_n^{2}\le n-n(1+\sqrt{\varepsilon})^{-1}\le n\sqrt{\varepsilon}.    
\end{equation}
 Furthermore, similar to Equation \eqref{BRg}, we obtain
\begin{equation} \label{BEc}
    \sum_{j\neq l}\frac{1}{\lambda_l\lambda_j}-R_g\le -2\sqrt\varepsilon (1+\sqrt{\varepsilon})^{-1}+\varepsilon^\frac{\alpha+1}{\alpha+2}\le -\sqrt\varepsilon \quad \text{on }\mathcal{B}\cap\mathcal{E}^c.
\end{equation}
Hence,
\begin{equation}\label{ab12}
    \begin{split}
        0\le& \int_{M^n} |\phi|^2\left[\sum_{j\neq l}\frac{1}{\lambda_j\lambda_l}-R_g\right]dV_g
        \\ \le & \int_{\mathcal{B}^c\cap \mathcal{E}^c}|\phi|^2  \varepsilon^\frac{\alpha+1}{\alpha+2}dV_g- \int_{\mathcal{B}\cap \mathcal{E}^c}|\phi|^2 \sqrt \varepsilon dV_g+\int_{\mathcal{E}}|\phi|^2 (n(n-1)-R_g)dV_g
        \\ \le &    \int_{\mathcal{B}^c\cap \mathcal{E}^c}|\phi|^2  \varepsilon^\frac{\alpha+1}{\alpha+2}dV_g- \int_{\mathcal{B}\cap \mathcal{E}^c}|\phi|^2  \sqrt\varepsilon dV_g+\left(\int_{\mathcal{E}}|\phi|^{\frac{2\alpha}{\alpha-2}}dV_g\right)^{\frac{\alpha-2}{\alpha}} \|(n(n-1)-R_g)_+\|_{\frac{\alpha}{2}}
        \\ \le & \varepsilon^\frac{\alpha+1}{\alpha+2} 2^6|M^n|- 
       \sqrt \varepsilon \int_{\mathcal{B}\cap \mathcal{E}^c}|\phi|^2  dV_g+
        \varepsilon (2^6)^\frac{\alpha-2}{\alpha}|M^n|.
    \end{split}
\end{equation}
Here, the first inequality follows from Llarull's integral formula, in the second inequality we used the definition of $\mathcal{E}$ and equation \eqref{BEc}, in the third inequality we applied H\"older inequality, and the fourth inequality follows from equation \eqref{C1al C2al}.

Combining H\"older inequality with equation \eqref{fal}, we obtain
\begin{equation} \label{abc}
    \begin{split}
        \int_{\mathcal{B}\cap \mathcal{E}^c}|\phi|^2 dV_g \ge& \int_{\mathcal{B}\cap \mathcal{E}^c}\left(\frac{3}{4}-3||\phi|-1|^2\right) dV_g
        \\ \ge & \frac{3}{4}|\mathcal{B}\cap \mathcal{E}^c|- 3\left(\int_{\mathcal{B}\cap \mathcal{E}^c}1 dV_g\right)^\frac{2}{\alpha} \left(\int_{\mathcal{B}\cap \mathcal{E}^c} ||\phi|-1|^\frac{2\alpha}{\alpha-2} dV_g\right)^{\frac{\alpha-2}{\alpha}}
        \\ \ge & \frac{3}{4}|\mathcal{B}\cap \mathcal{E}^c|-
       3|\mathcal{B}\cap \mathcal{E}^c|^\frac{2}{\alpha}\left( \varepsilon_*^\frac{\alpha}{\alpha-2}|M^n|\right)^{\frac{\alpha-2}{\alpha}}
       \\ \ge & \frac{1}{2}|\mathcal{B}\cap \mathcal{E}^c|- \frac{\alpha-2}{\alpha} \left[3^\frac{\alpha}{\alpha-2}\left(\frac{8}{\alpha}\right)^\frac{2}{\alpha-2}\right]\varepsilon_*^\frac{\alpha}{\alpha-2}|M^n|
       \\ \ge & \frac{1}{2}|\mathcal{B}\cap \mathcal{E}^c|-192\varepsilon_*^\frac{\alpha}{\alpha-2}|M^n|
    \end{split}
\end{equation}
where the second last inequality follows from Young's inequality which states that for any $x,y>0$,
\begin{equation}
    \frac{2}{\alpha}\left(\frac{\alpha}{8}x\right)+\frac{\alpha-2}{\alpha}\left[3^\frac{\alpha}{\alpha-2}\left(\frac{8}{\alpha}\right)^\frac{2}{\alpha-2}y\right]
    \ge \left(\frac{\alpha}{8}x\right)^\frac{2}{\alpha}\left[3^\frac{\alpha}{\alpha-2}\left(\frac{8}{\alpha}\right)^\frac{2}{\alpha-2}y\right]^\frac{\alpha-2}{\alpha}=3x^\frac{2}{\alpha}y^\frac{\alpha-2}{\alpha}.
\end{equation}
Combining Equation \eqref{ab12} and \eqref{abc}, we have
\begin{equation}
    |\mathcal{B}\cap \mathcal{E}^c|\le \left[128\varepsilon^\frac{\alpha}{2(\alpha+2)}+384\varepsilon_*^\frac{\alpha}{\alpha-2}+\varepsilon^\frac{1}{2}2\cdot (64)^\frac{\alpha-2}{\alpha}\right]|M^n|\le \left[256\varepsilon^\frac{\alpha}{2(\alpha+2)}+384\varepsilon_*^\frac{\alpha}{\alpha-2}\right]|M^n|.
\end{equation}
With the help of \eqref{B estimate}, we obtain
\begin{equation} \label{Acg}
    |\mathcal{B}|\le |\mathcal{E}|+|\mathcal{B}\cap \mathcal{E}^c|\le \left[257\varepsilon^\frac{\alpha}{2(\alpha+2)}+384\varepsilon_*^\frac{\alpha}{\alpha-2}\right]|M^n|
\end{equation}
which finishes the proof. 
Moreover,  we have $|\mathcal{B}^c|_g\le  (1+\sqrt\varepsilon)^\frac{n}{4}|S^n|_{g_0}$.
Hence,
\begin{equation} \label{volume converge}
    |S^n|_g\le \left\{1-\left[257\varepsilon^\frac{\alpha}{2(\alpha+2)}+384\varepsilon_*^\frac{\alpha}{\alpha-2}\right]\right\}^{-1}(1+\sqrt\varepsilon)^\frac{n}{4}|S^n|_{g_0}.
\end{equation}
The estimate above will be used in the proof of Theorem \ref{thm main 2} below.
\end{proof}

\begin{remark} \label{example j}
Given a family of manifolds $(M_i,g_i)$ satisfying the assumptions of Theorem \ref{thm main 3} with $\varepsilon=\tfrac1i$, we obtain $C^0$ subsequential convergence $g_i\to g_0$ outside an arbitrary small set.
We remark that taking a subsequence is necessary, and it is easy to construct counter examples otherwise.
\end{remark}

\begin{proof}[Proof of Theorem \ref{thm main 2}] 
We follow the argument of B.~Allen, E.~Bryden and D.~Kazaras \cite{ABK}.
    In view of \cite[Theorem 1.1]{APS} by B.~Allen, R.~Perales and C.~Sormani, it suffices to show volume convergence.
    By assumption, $g_i\ge g_0$ which implies $|S^n|_{g_i}\ge|S^n|_{g_0}$.
      Combining this with estimate \eqref{volume converge}, we obtain
$|S^n|_{g_i}\to |S^n|_{g_0}$ which finishes the proof.
\end{proof}

\appendix

\section{Llarull's integral formula}\label{appendix:llarull}

\begin{theorem} \label{A.1}
Let $(M^n,g)$ be a smooth compact spin manifold without boundary. 
Let $f:M^n\to S^n$ be a smooth area non-increasing map.
Consider a harmonic spinor $\psi$ in the twisted spinor bundle $\mathcal S(M^n)\otimes f^\ast \mathcal S(S^n)$.
Then
\begin{equation}
    \int_{M^n} |\phi|^2\left[\sum_{j\neq l}\frac{1}{\lambda_j\lambda_l}-R_g\right]dV_g\ge 4\int_{M^n} |\nabla \phi|^2dV_g
\end{equation}
where $\lambda_i$ are the square roots of eigenvalues of $g$ with respect to $f^\ast g_0$.
Moreover, if $n$ is an even integer, the index of the Dirac operator $\slashed D$ is non-vanishing, and we can always find a non-trivial $\phi \in \mathcal S(M^n)\otimes f^\ast \mathcal S(S^n)$ solving $\slashed D\phi=0$.
\end{theorem}

Such a formula also holds for target manifolds other than the sphere as long as they have positive curvature operator and non-vanishing Euler characteristic, cf. \cite{GoetteSemmelmann} by S.~Goette and U.~Semmelmann.
We expect similar stability results to hold in this setting.

\subsection{Odd dimensional case}\label{SS:odd}
Theorem \ref{A.1} has been the crucial tool to prove Theorem \ref{New main theorem} and Theorem \ref{thm main 3} in even dimensions.
We now describe the required adjustments to also obtain the odd dimensional case.
This closely follows Llarull's original ideas \cite{Llarull}.

Let $n\ge3$ be an odd integer.
Consider the product metric $\overline{g}= ds^2+g$ on $M^{n}\times S_r^1$, where $S_r^1$ is a circle of radius $r$.
We note that $R_g=R_{\overline{g}}$ and 
\begin{equation}
   M^n\times S^1_r\xrightarrow{f\times\operatorname{Id}} S^n\times S_r^1\xrightarrow{\operatorname{Id}\times \frac{1}{r}\operatorname{Id}}S^{n}\times S^1\xrightarrow{h} S^n\wedge S^1 \cong S^{n+1},  
\end{equation}
where $\wedge$ is the smash product and $h$ is a 1-contracting map.
Applying the Llarull's formula to the even dimensional manifold $(M^n\times S^1_r,\overline{g})$ and using Kato's inequality, we obtain
\begin{equation}
    \int_{M^n\times S^1_r}\left[2\sum_{1\le l<j\le n}\frac{1}{\lambda_l\lambda_j}+\frac{2}{r}\sum_{l=1}^n\frac{1}{\lambda_l}-R_g\right]|\phi|^2_{
   \overline{g}
    }dV_{\overline{g}}\ge \int_{M^n\times S^1_r} 4|\nabla \phi|^2_{\overline{g}}dV_{\overline{g}}\ge \int_{M^n\times S^1_r} 4|\nabla |\phi|_{\overline{g}}|^2_{\overline{g}}dV_{\overline{g}}.
\end{equation}
For any $\delta>0$, by choosing $r$ sufficiently large, there exists a point $p$ on  $S^1_r$ such that
\begin{equation}\label{odd}
    \int_{M^n\times \{p\}}\left[2\sum_{1\le l<j\le n}\frac{1}{\lambda_l\lambda_j}-R_g-\delta\right]|\phi|_{\overline{g}}^2dV_{g}\ge \int_{M^n\times \{p\}} 4|\nabla |\phi|_{\overline{g}}|^2_{g}dV_{g}.
\end{equation}
Hence, we may continue as before and deduce that Theorem \ref{New main theorem} and Theorem \ref{thm main 3} also holds in odd dimensions.

\end{document}